\documentclass[10pt]{amsart}

\usepackage[margin=34mm, marginpar=25mm]{geometry}
\usepackage[utf8]{inputenc}
\usepackage{tikz-cd}
\usepackage{amsfonts, amstext, amsmath, amsthm, amscd, amssymb}
\usepackage{mathrsfs}
\usepackage{tikz}
\usetikzlibrary{knots}
\tikzset{point/.style={insert path={ node[scale=2.5*sqrt(\pgflinewidth)]{.} }}}
\usepackage{url}
\usepackage[all]{xypic}

\usepackage{color}
\usepackage[colorlinks,
    linkcolor={red!50!black},
    citecolor={blue!50!black},
    urlcolor={blue!80!black}]{hyperref}

\usepackage{enumerate}

\newcommand{\bp}{\begin{pmatrix}}
\newcommand{\ep}{\end{pmatrix}}
\newcommand{\be}{\begin{equation}}
\newcommand{\ee}{\end{equation}}

\newcommand{\smfrac}[2]{\mbox{\footnotesize$\displaystyle\frac{#1}{#2}$}} % small medium frac
\numberwithin{equation}{section}

\theoremstyle{plain}

\newtheorem*{theorem*}{Theorem}

\theoremstyle{definition}

\theoremstyle{remark}

\newtheorem*{remark*}{Remark}

\numberwithin{equation}{section}

\def\Z{\mathbb Z}

\def\Q{\mathbb Q}

\def\bp{\begin{pmatrix}}
\def\ep{\end{pmatrix}}
\def\ba{\begin{array}}
\def\ea{\end{array}}
\def\bn{\begin{enumerate}}
\def\en{\end{enumerate}}

\DeclareMathOperator\Wh{Wh}

\DeclareMathOperator\Bl{Bl}

\begin{document}

\title[Non-approximable slice discs]{Smoothly slice knots with smoothly non-approximable topological slice discs}

\author{Min Hoon Kim}
\address{Department of Mathematics, Ewha Womans University, Seoul, Republic of Korea}
\email{minhoonkim@ewha.ac.kr}
\author{Mark Powell}
\address{School of Mathematics and Statistics,
University of Glasgow, University Place, Glasgow, G12 8QQ, United Kingdom}
\email{mark.powell@glasgow.ac.uk}

\def\subjclassname{\textup{2020} Mathematics Subject Classification}
\expandafter\let\csname subjclassname@1991\endcsname=\subjclassname
\expandafter\let\csname subjclassname@2000\endcsname=\subjclassname

\begin{abstract}
We construct infinitely many smoothly slice knots having topological slice discs that are non-approximable by smooth slice discs.
\end{abstract}

\maketitle

Let $K$ be a smoothly slice knot in $S^3$ and let $\mathcal{D}$ be a topological slice disc for $K$ in $D^4$, that is, $\mathcal{D}$ is a locally flat, topologically embedded disc in $D^4$ such that $\partial (D^4, \mathcal{D})=(S^3,K)$.  We say that $\mathcal{D}$ is \emph{smoothly non-approximable} if there exists $\varepsilon>0$ such that for every topological slice disc $D$ that is topologically isotopic rel.\ boundary to $\mathcal{D}$, the $\varepsilon$-neighbourhood $\nu_{\varepsilon}(D)$  of $D$ does not contain any smooth slice disc for $K$.

\begin{theorem*}\label{theorem:main}
There exist infinitely many smoothly slice knots $K_m$ in $S^3$, each of which admits a smoothly non-approximable topological slice disc.
\end{theorem*}

This answers a refined version of mathoverflow question \url{https://mathoverflow.net/questions/419693}, asked by M.~Winter.

It is interesting to contrast our theorem with Venema's results in~\cite{Venema-approximating-disks,Venema-approximating}, which imply that every topological slice disc can be approximated by smoothly embedded disc if the boundary is allowed to move. In particular Venema's smooth approximation could have a different knot type on the boundary, or could even lie in the interior of $D^4$.

\begin{figure}[htb!]
\centering
\includegraphics[width=0.8\textwidth]{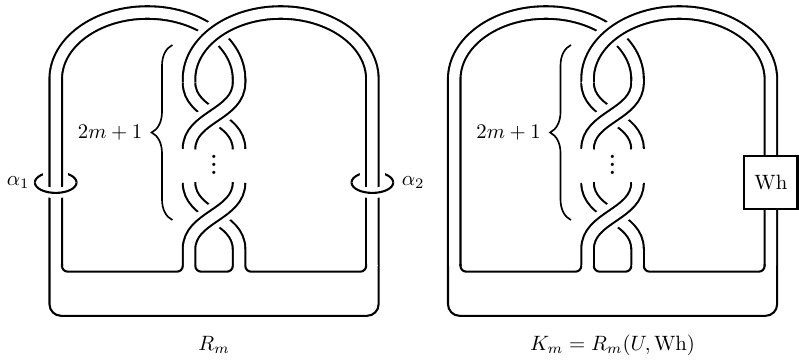}
\caption{The knots $R_m$ and $K_m$. There are $2m+1$ crossings between the bands. On the left, the knot $R_m$ has unknots $\alpha_1$ and $\alpha_2$ in its complement, each linking one of the bands of the obvious Seifert surface for $R_m$. Perform a satellite operation on $\alpha_2$ using $\Wh$ to obtain $K_m := R_m(U,{\Wh})$, as shown on the right. } \label{figure1}
\end{figure}

\begin{proof}
For an odd integer $m\geq 1$, let ${\Wh} := \operatorname{Wh}(T_{2,3})$ denote the positive-clasped zero-twisted Whitehead double of the right-handed trefoil, let $U$ denote the unknot, and let $K_m = R_m (U, {\Wh})$, as defined in the caption of Figure~\ref{figure1}. Since $\alpha_1$ and $\alpha_2$ link $R_m$ trivially, the Alexander polynomial is unchanged by the satellites~\cite{Seifert:1950-1}, i.e.\ $\Delta_{K_m}(t) = \Delta_{R_m}(t)$.  Also $\Delta_{R_m}(t) = ((m+1)t - m)(mt-(m+1))$. Thus if $m\neq n$, then $\Delta_{K_m}(t) \neq \Delta_{K_n}(t)$, and hence the knots $K_m$ are mutually distinct.

The configuration $(R_m,\alpha_1,\alpha_2)$ is often considered as an operator $R_m(-,-)$ that takes two knots as input, with $R_m(J_1,J_2)$ obtained by performing the satellite operation on $\alpha_i$ with $J_i$, for $i=1,2$. Since taking $J_1=U$ has no effect, this explains our definition in the caption to Figure~\ref{figure1}.

For each $m$, the knot~$K_m$ is smoothly slice, as can be seen by `cutting' the right hand band in the Seifert surface. From the embedded Morse theory perspective, this corresponds to a saddle move, which results in a 2-component unlink, which can then be capped off by two smooth discs, corresponding to minima, to obtain a smooth slice disc for $K_m$. \label{smooth-slice-disc}

Let $\mathcal{D}_m$ be a topological slice disc for $K_m$ obtained by cutting the left band and capping with two parallel copies of a topological slice disc for ${\Wh}$, which exists since $\Delta_{\Wh}(t) =1$~\cite{Freedman:1984-1}.
We claim that $\mathcal{D}_m$ is smoothly non-approximable for each $m$.

For a contradiction, supposed that for every $\varepsilon>0$, there exists a topological slice disc $D_m$  for $K$ that is topologically isotopic rel.\ boundary to $\mathcal{D}_m$, such that the $\varepsilon$-neighbourhood $\nu_{\varepsilon}(D_m)$ contains a smooth slice disc $D_{m,\varepsilon}'$ for $K$.
Let $N(D_m)$ be an open tubular neighbourhood of the topologically locally flat slice disc $D_m$, which exists by~\cite[Theorem~9.3]{Freedman-Quinn:1990-1}. Then $D^4 \smallsetminus N(D_m)$ is the exterior of $D_m$, and $\partial(D^4 \smallsetminus N(D_m)) \cong M_{K_m}$, the 3-manifold obtained by zero-framed surgery on $K_m$.
Choose $\varepsilon>0$ small enough such that $\nu_\varepsilon(D_m) \subseteq N(D_m)$.
We have inclusion maps
\[M_{K_m} \xrightarrow{i_1} D^4\smallsetminus N(D_m) \xrightarrow{i_2} D^4\smallsetminus \nu_{\varepsilon}(D_m) \xrightarrow{i_3} D^4\smallsetminus D_{m,\varepsilon}'.\]
It follows from Alexander duality that these inclusion maps induce isomorphisms
\[H_1(M_{K_m};\Z) \xrightarrow{\cong} H_1(D^4 \smallsetminus N(D_m);\Z) \xrightarrow{\cong} H_1(D^4\smallsetminus D_{m,\varepsilon}';\Z).\]
Since $H_1(M_{K_m};\Z) \cong \Z$ we obtain consistent coefficient systems $\pi_1(X) \to \Z \cong \langle t \rangle$ for all $X \in \{M_{K_m}, D^4 \smallsetminus N(D_m),D^4\smallsetminus D_{m,\varepsilon}'\}$, and can therefore consider homology with $\Q[t^{\pm 1}]$ coefficients.
Consider the corresponding composition
\[\iota \colon H_1(M_{K_m};\Q[t^{\pm 1}]) \xrightarrow{(i_1)_*} H_1(D^4\smallsetminus N(D_m);\Q[t^{\pm 1}]) \xrightarrow{(i_3 \circ i_2)_*} H_1(D^4\smallsetminus D_{m,\varepsilon}';\Q[t^{\pm 1}]).\]
We need to understand the kernel of this map $\iota := (i_3\circ i_2\circ i_1)_*$.

By our description of $D_m$, $(i_1)_*(\alpha_1)=0$ and hence $\langle \alpha_1 \rangle \subseteq \ker \iota$, where we abuse notation and use $\alpha_1$ to also denote the corresponding curve in the complement of~$K_m$. Since $D_{m,\varepsilon}'$ is a slice disc, the kernel of $\iota$ is a metaboliser for the Blanchfield form of~$M_{K_m}$ (see e.g.~\cite[Theorem 4.4]{Cochran-Orr-Teichner:1999-1} or \cite[Theorem~2.4]{Hillman:2012-1-second-ed}).
The Blanchfield form is a Hermitian, sesquilinear, nonsingular pairing
\[\Bl \colon H_1(M_{K_m};\Q[t^{\pm 1}]) \times H_1(M_{K_m};\Q[t^{\pm 1}]) \to \Q(t)/\Q[t^{\pm 1}],\]
and we say that $P \subseteq H_1(M_{K_m};\Q[t^{\pm 1}])$ is a metaboliser if $P = P^{\bot}$.
We can compute, using the presentation matrix $tV_m -V_m^T$, where $V_m = \left[\begin{smallmatrix} 0 & m+1 \\ m & 0 \end{smallmatrix}\right]$ is a Seifert matrix for~$K_m$, that  \[H_1(M_{K_m};\Q[t^{\pm 1}]) \cong  \smfrac{\Q[t^{\pm 1}]}{(mt-(m+1))}  \oplus \smfrac{\Q[t^{\pm 1}]}{((m+1)t - m)} \cong \Q \oplus \Q.\]
The first summand is generated by $\alpha_1$, and the other is generated by $\alpha_2$.

We claim that $\langle \alpha_1 \rangle=\ker(\iota)$. To see this, we consider $H_1(M_{K_m};\Q[t^{\pm 1}])$ as a $\Q$-vector space. Then $\langle \alpha_1\rangle\subseteq \ker(\iota) \subseteq H_1(M_{K_m};\Q[t^{\pm 1}])$ implies that $1 \leq \dim_\Q \ker(\iota) \leq 2$.
However $\ker(\iota)$ is a metaboliser for the Blanchfield form, and the Blanchfield form is nonsingular, so we cannot have $\ker \iota = H_1(M_{K_m};\Q[t^{\pm 1}])$. Hence $\dim_\Q \ker(\iota)=1$. Since $\langle \alpha_1\rangle\subseteq \ker(\iota)$ and their dimensions as $\Q$-vector spaces are equal, we deduce that $\langle \alpha_1\rangle = \ker(\iota)$ as claimed.

Now we apply the $d$-invariant from Heegaard-Floer homology, together with calculations of Cha-Kim~\cite{Cha-Kim:2021-1} and Cha~\cite{Cha:2021-1}, to deduce that the smooth disc $D_{m,\varepsilon}'$ cannot exist.
For the definition of the $d$-invariant we refer to Ozsv\'ath--Szab\'o~\cite[Definition~4.1]{Ozsvath-Szabo:2003-2}.  It suffices for our purposes to know that given a rational homology 3-sphere $Y$ and a spin${}^c$ structure $\mathfrak{s}$ on $Y$, the $d$-invariant is a rational number $d(Y,\mathfrak{s}) \in \Q$.

Let $\Sigma_r$ be the $r$-fold cyclic cover of $S^3$ branched along $K_m$, which is a rational homology 3-sphere for all prime $r \geq 2$. Then $\alpha_1$ and $\alpha_2$ lift to homology classes in $H_1(\Sigma_r;\Z)$, say $x_1$ and $x_2$ respectively. Since $H_1(\Sigma_r;\Z)\cong H_1(M_{K_m};\Z[t^{\pm 1}])/(t^r-1)$, we have that
\[H_1(\Sigma_r;\Z)\cong \Z/_{(m+1)^r-m^r}\oplus  \Z/_{(m+1)^r-m^r},\] where the summands are generated by $x_1$ and $x_2$ respectively. In particular, $\Sigma_r$ has a unique spin structure $\mathfrak{s}_{\Sigma_r}$ since it is a $\Z/2$-homology 3-sphere.

By \cite[Lemma~5.2]{Cha-Kim:2021-1} and \cite[p.\ 17, Assertion]{Cha:2021-1}, the condition $\langle \alpha_1 \rangle=\ker(\iota)$ implies that there exists a prime $r$ such that the kernel of the inclusion induced map $H_1(\Sigma_{r};\Z)\to H_1(V_r;\Z)$ is generated by $x_1$, where $V_r$ is the $r$-fold cyclic cover of $D^4$ branched along $D'_{m,\varepsilon}$. Since $D'_{m,\varepsilon}$ is a smooth slice disc, \cite[Theorem 1.1]{Grigsby-Ruberman-Strle:2008-1} implies that
$$d(\Sigma_r, \mathfrak{s}_{\Sigma_r}+k\widehat{x}_1)=0$$
 for all $k\in \Z$.  Here the spin structure $\mathfrak{s}_{\Sigma_r}$ uniquely determines a spin${}^c$ structure, and $\widehat{x}_1 := PD^{-1}(x_1) \in H^2(\Sigma_r;\Z)$ is the Poincar\'{e} dual of $x_1$. Then the spin${}^c$ structure $\mathfrak{s}_{\Sigma_r}+k\widehat{x}_1$ is defined by noting that spin${}^c$ structures on $\Sigma_r$ are a torsor over $H^2(\Sigma_r;\Z)$.
% where $\mathfrak{s}_{\Sigma_r}$ is the spin structure of $\Sigma_r$.

  We obtain the desired contradiction since there exists an integer $k\in \Z$ such that $$d(\Sigma_r, \mathfrak{s}_{\Sigma_r}+k\widehat{x}_1)\neq 0,$$ by \cite[Theorem~5.4]{Cha-Kim:2021-1} and \cite[Lemma 4.1]{Cha:2021-1}. This completes the proof.
\end{proof}

\begin{remark*}
  Our examples are topologically doubly slice and smoothly slice, but not smoothly doubly slice. Infinitely many such knots were first constructed by Meier~\cite{Meier-smooth-vs-top-doubly-slice}, though the earlier $d$-invariant arguments of Cochran, Harvey, and Horn \cite[pp.\ 2140--1]{Cochran-Harvey-Horn:2013-1} implicitly show the existence of such a knot.  These works do not suffice to prove our theorem, because for a doubly slice knot $K$ each individual slice disc  has the property that for each $r$,  the map $H_1(\Sigma_r;\Z) \to H_1(V_r;\Z)$ is surjective. Moreover the kernels of these maps, for the  two $V_r$ corresponding to each slice disc, are two complementary summands of $H_1(\Sigma_r;\Z)$. Then, as in \cite{Cochran-Harvey-Horn:2013-1,Meier-smooth-vs-top-doubly-slice}, one can restrict to a single prime $r$, and compute the invariants $d(\Sigma_r,\mathfrak{s})$ to obtain a contradiction. In our case, we cannot assume any `homology ribbon' property, so we cannot apply the same argument. The putative slice disc~$D_{m,\varepsilon}'$ and the actual smooth slice disc for $K_m$ from the top of page~\pageref{smooth-slice-disc} could a priori share the same kernel. We apply the results of  \cite{Cha-Kim:2021-1,Cha:2021-1}, which circumvent this subtle issue by considering $d$-invariants of~$\Sigma_r$ for arbitrarily large primes $r$. Similar $d$-invariant arguments were applied in Kim-Livingston~\cite{Kim-Livingston:2022-1} to construct an infinite family of topologically slice knots that are not smoothly concordant to their reverses.
\end{remark*}

\begin{remark*}
The fact that an isotopy from $\mathcal{D}$ to $D$ is permitted in the definition of smoothly non-approximable distinguishes our examples from topological slice discs constructed by the following argument, told to us by Robert Gompf.
Start with a smooth slice disc $\Delta$ for a knot $K$.  For any $\delta>0$, in the $\delta$-neighbourhood of $\Delta$ we may construct a capped tower with one storey $T$, as in \cite{Freedman-Quinn:1990-1}, with attaching circle $K$, such that $T$ does not contain any smooth disc with boundary~$K$. Such a $T$ exists, as shown in e.g.~\cite{Gompf:1984-1}.  For some $\varepsilon < \delta$, there is a topological slice disc $D$ for $K$ and $\nu_\varepsilon(D)$ contained within $T$, with the same attaching circle~\cite{Freedman:1982-1,Freedman-Quinn:1990-1}.  There cannot be any smooth slice disc for $K$ within $\nu_\varepsilon(D)$, or else there would also be such a smooth slice disc lying within $T$. This provides an alternative answer to Winter's original question, which did not allow an isotopy.

However, it can be shown that for every $\varepsilon>0$,  $D$ is topologically isotopic to a disc that contains $\Delta$ in its $\varepsilon$-neighbourhood, and hence $D$ fails to be smoothly non-approximable in the sense of our definition.
\end{remark*}

\subsubsection*{Acknowledgements}
We are grateful to Martin Winter for asking the question on mathoverflow that led to this note, and to both Chuck Livingston and an anonymous referee for helpful comments. We are also grateful to Bob Gompf for pointing out an alternative answer to Winter's question discussed in the final remark, which led us to allow an isotopy in the definition of smoothly non-approximable.
Min Hoon Kim was partially supported by the Samsung Science and Technology Foundation (SSTF-BA2202-01) and National Research Foundation grant 2021R1C1C1012939.
Mark Powell was partially supported by EPSRC New Investigator grant EP/T028335/2.

\bibliographystyle{amsalpha}
\def\MR#1{}
\bibliography{research}
\end{document}